# A New Decision- Making Method Based on Shannon Entropy Analysis

H. Babaei[1*], SH. Mohammadi[2], H. Ghaneai[3]

[1] *Department of Management, Meybod University, Meybod, Iran*

[2] *Department of accounting, economic and management, Yazd University, Yazd, Iran*

[3] *Department of Computer Engineering, Meybod University, Meybod, Iran*

**Abstract**

Because the appropriate combination of existing elements and establishing coordination between them as a consequence of making the right decision to accomplish the intended objective is achieved, management is now one of the main pillars of community management. Decisions are made in the majority of situations when the decision maker is pleased with the conclusion based on numerous factors. Several criteria are used instead of one measure of optimality in multi-criteria decision making, which has been studied by numerous academics in recent decades. The importance of the indicators in this type of decision making is clearly not equal, and it is necessary to understand the coefficient of importance or weight of each of these indicators in decision making. In this work, a novel technique termed scattering axis Dispersion-based Weighting Method (**DWM**) is suggested to address weighing issues, with the closest method in terms of computational logic being their entropy. After constructing the option's criterion matrix, the mean, standard deviation, and coefficient of variation are determined, and then the weight of each criterion is calculated, according to the proposed **DWM** technique. Several numerical examples have been utilized to demonstrate and assess the suggested technique. In addition, the Shannon entropy approach, which is a commonly used weighting method, was chosen to compare the findings. The statistical results demonstrate that these two weighing techniques have a strong connection. In compared to the Shannon entropy technique, the suggested method has the following advantages:

- Minimal computational burden
- The data does not need to be normalized.
- Its use in the case of negative data.

*Corresponding author: Tel: +, E-mail: babaei@meybod.ac.ir

# 1. Introduction

Choosing and making the ideal option has always been a difficult task. When making a decision is difficult due to factors such as the availability of several options and the examination of numerous indications, the use of Multiple Criteria Decision Making (MCDM) procedures is recommended [1], which is a highly significant area of decision-making theory [2]. Several criteria are employed instead of one measure of optimality in multi-criteria decision making approaches that have been studied by scholars in recent decades. The phases of the problem-solving and decision-making process will lead to each organization's direction toward attaining its objectives. The organization's management occupies a crucial position in the market. Because the quality of plans and programs, the efficacy and efficiency of strategies, and the quality of the outcomes gained from their implementation all depend on the quality of the manager's decisions, management quality is basically a function of decision-making quality. When a choice is made based on many factors, it is usually desirable and satisfying to the decision maker [1]. Having several stakeholders with varying interests and values further complicates the decision-making process on various choices [3].

Indicators, on the other hand, are not all created equal when it comes to making decisions. In such instances, the significance of these indicators must be determined, and the coefficient of relevance or weight of each of these indicators in decision-making must be determined. Each index's weight reflects how important it is in comparison to other indicators [4]. The deliberate and accurate selection of weights will go a long way toward obtaining the intended result. Although a significant portion of MCDM approaches were developed in the previous century, for example [5-7], their widespread usage and occasional flaws drive scientists to design new ways or alter current ones, as shown in [8-10].

There are two types of multi-criteria decision making (MCDM) methods: Multi Objective Decision Making (MODM) [3] and Multi Attribute Decision Making (MADM) [11]. The goal of decision making is to select the best choice or to balance the many decision variables. Each technique of decision-making has an unique task: one is to weight the criteria, another is to rank the alternatives, and the third is to assess the criteria. Multi-objective models are used for design in general, and multi-objective models are utilized to choose the superior choice, which is the major emphasis of this article. The primary distinction between multi-objective and multi-criteria decision models is that the former is specified in a continuous decision space while the later is defined in a discrete decision space [2].

The relevance of indicators in decision-making, on the other hand, is clearly not equal. In such instances, the significance of these indicators must be determined, and the coefficient of relevance or weight of each of these indicators in decision-making must be determined. Each index's weight reflects how important it is in comparison to other indicators [5]. The deliberate and accurate selection of weights will go a long way toward obtaining the intended result. The Analytic Hierarchy Process technique (AHP) [6, 7] is one of several ways for calculating the weight of criteria. The network analysis procedure approach (ANP) [8], the entropy weighting method (EWM) [9-11], and the Criteria Importance Through Intercriteria Correlation method (CRITIC) [12] and the Bulls-eye method [13] are all commonly utilized methods. We can employ the Superiority and Inferiority Ranking (SIR) [14], the Step wise Weight Assessment Ratio Analysis (SWARA) [15], the Best Worst method (BWM) [2], and the Multi-attribute evaluation using importance weight estimates method (IMP) [16] as some of the novel decision-making methods. Although a significant portion of MCDM approaches were created in the previous century, for example [6, 7, 17], their widespread usage and occasional flaws drive scientists to propose new ways or alter current ones, for example [18-20].

One of the most commonly used methods for ranking and evaluating the relevance of elements that prioritize each criteria using pairwise comparisons of possibilities is the hierarchical analysis approach [21]. This method is also employed for weighing criteria. Because paired comparisons get more complex as the number of items in each cluster grows, decision criteria are generally split into sub-criteria. Zolésio [22] introduced the CRITIC method for determining the weight of the criterion. The opinion of specialists is less relevant in this technique. One of the crucial method's most essential capabilities is that specialists are not involved in the process. The data is evaluated using this approach depending on the level of interference and conflict between factors or criteria [23]. The best-worst method (BWM) is a multi-criteria decision-making approach that is among the new multi-criteria decision-making strategies [24]. The decision maker determines the best and worst indicators, and a pair comparison is conducted between each of these two indicators (best and worst) and other indicators in this approach. Then, in order to estimate the weight of various indicators, a maximal problem is constructed and solved. In order to evaluate the correctness of the comparisons, a formula for determining the incompatibility rate was also explored in this technique [2].

One of the multi-criteria decision-making approaches for determining weight in a three-parameter gray number spectrum is the Bulls-eye methodology. In decision matrices, three parameters are utilized for weighting. With little usage of experts and recalculation of weights, the Bulls-eye approach may substantially decrease the impact of human judgements on criterion weighting.

One of the multi-criteria decision-making approaches, the Swara method or progressive weight evaluation ratio analysis, tries to determine the weight of criterion and sub-criteria [25]. The criteria are ordered according to their worth in this way. The most significant criterion is listed first, and the least important criterion is ranked last in this technique. Experts (respondents) play an essential role in establishing the weight of the criteria in this technique. The ability to estimate experts and experts in relation to the relevance of criteria in the process of establishing their weight is the method's key characteristic. This approach may be used to collect and organize information from experts and intellectuals [26]. In information theory, entropy is used to quantify a message's information content [27]. Entropy may be used to assess characteristics since the decision matrix for a set of alternatives includes a certain amount of information [28]. The entropy technique is a multi-criteria decision-making approach based on dispersion that determines the weight of the criteria depending on a range of criteria [29]. The basic premise of this technique is that the greater the spread in an index's values, the more significant that index is [29]. We can remove that criteria when all alternatives have the same values for that criterion [30].

The following matrix form can be used to completely characterize a multi attribute decision making:

$$B = \begin{array}{c} \\ a_1 \\ a_2 \\ \vdots \\ a_m \end{array} \begin{array}{c} c_1 \quad c_2 \quad \cdots \quad c_n \\ \left( \begin{array}{cccc} p_{11} & p_{12} & \cdots & p_{1n} \\ p_{21} & p_{22} & \cdots & p_{2n} \\ \vdots & \vdots & \ddots & \vdots \\ p_{m1} & p_{m2} & \cdots & p_{mn} \end{array} \right) \end{array}$$

The ith option is represented by $a_i$, the jth index is represented by $C_j$ and the value of the ith option in terms of the jth index is represented by $p_{ij}$. The objective is to pick the best choice (the one that is most desirable and important), or an alternative that has the best overall value. Various approaches may be used to determine the alternative value of alternative $v_i$. In general, if we give the weight wj to the criteria j ($w_j \geq 0, \sum w_j = 1$)، we may get $v_i$ using the simple weighting technique [31], which is the foundation for many MCDM systems, as follows:

$$v_i = \sum_{j=1}^{n} w_j \, p_{ij}$$

The approach for obtaining standard or vector weights $w = \{w_1, w_2, \ldots, w_n\}$ [2] is highly essential here, and it is the rationale for creating various MCDM systems in recent decades.

The deliberate and accurate selection of weights will go a long way toward obtaining the intended result. The entropy technique may be used to calculate the weight of the indicators once the data for a decision matrix has been properly provided. In the social sciences, physics, and information theory, entropy is a crucial term. A criterion-option matrix is required for this approach. Shannon and Weaver suggested this technique in 1974. The degree of uncertainty in a continuous probability distribution is represented by entropy. The basic premise behind this technique is that the greater the dispersion in an index's values, the more important that index is. Shannon demonstrated that events with a high likelihood of occurrence offer less information, whereas occurrences with a low chance of occurrence provide more information. The uncertainties are really decreased when new information is acquired, and the value of the new information is equal to the amount of uncertainty that is reduced. As a result, uncertainty and information are two factors that are linked. In many situations, the weights determined using their entropy technique appear unreasonable. There are other steps in this approach, as well as a large number of computations.

The goal of this work is to offer a novel weighting technique (WBI) that is the most computational logic entropy efficient. We shall explain that one of the advantages of the weighting technique based on this logic is its low computing burden when compared to the Shannon entropy method, as well as the absence of the requirement to normalize the data. As a result, understanding this approach can assist researchers in resolving choice issues involving a huge number of variables. The following is the continuation of this article. A novel technique (**DWM**) is proposed in Section 2. In Section 3, **DWM** is applied to a real-world situation and its entropy is compared to it using different assessment criteria. Section 4 contains the conclusions and recommendations for future study.

## 2. Dispersion-based Weighting Method (DWM)

The weight of the components is critical in most multi-criteria decision making. Entropy is a concept used in social sciences, physics, and information theory, as well as one of the multi-criteria decision-making approaches for determining the weight of factors. This method necessitates the creation of a matrix based on criteria and choices. The entropy approach may be used to assess the weights if the decision matrix data is known. According to the weights acquired from the indices, those indices with higher dispersion are more significant than other indices in the Shannon entropy approach, and their influence is greater in picking the best choice [32].

**2.1. Shannon entropy method**

**Step 1:** The decision matrix is created first. To create this matrix, just determine whether the criteria are qualitative, assess each choice in respect to each criterion using verbal phrases, and, if the criteria are small, enter the real number of that evaluation. The criteria are listed in the columns, and the alternatives are shown in the rows in the decision matrix shown below. The score of the first choice, for example, is $x_{12}$ when compared to the second criterion.

$$X_{ij} = [x_{ij}]_{n \times m} = \begin{bmatrix} x_{11} & x_{12} & \cdots & x_{1m} \\ x_{21} & x_{22} & \cdots & x_{2m} \\ \vdots & \vdots & \vdots & \vdots \\ x_{n1} & x_{n2} & \cdots & x_{nm} \end{bmatrix} \quad (1)$$

**Step 2:** Each normalized flow is called after the aforementioned matrix is normalized. Normalization is done linearly in the Shannon entropy approach. In this method, each column's knowledge is split by the total of that column's knowledge.

$$[P_{ij}]_{n \times m} = \left[ x_{ij} / \sum_{i=1}^{n} x_{ij} \right]_{n \times m} \quad (\text{٢})$$

**Step 3:** Calculating the entropy of each index: Shannon used the following formula to calculate the entropy of a probability distribution for each random phenomenon:

$$E = S\begin{pmatrix} P_1 \\ P_2 \\ \vdots \\ P_m \end{pmatrix}, \sum_{i=1}^{m} P_i = 1$$

The following formula was proposed to compute the entropy of such events, which, due to the uncertainty of the numbers within the matrix, also contain the indice0s:

$$E_j = -k \sum_{i=1}^{m} [P_i . \ln P_i], k = \frac{1}{\ln(m)} \quad (3)$$

Where Ej is the index's entropy is j, the number of choices is m. From the perspective of the i, Pi option, the likely value of the index value, ln is the symbol for the neper logarithm or natural logarithm, and k is a constant number for modifying entropy between 0 and 1.

It's worth noting that in decision matrices, it's usually m ≥3, which means that if there are fewer than three alternatives, it's not as essential, so:

$$\frac{1}{\ln(m)}, (m = 3 > e = 2.7 \rightarrow \ln(m) > 1 \rightarrow \frac{1}{\ln(m)} < 1)$$

In this formula, the closer Ej, the jth index's entropy, goes to one, the less the index's influence on prioritizing alternatives is decreased to zero. As a result, if a phenomena or index has the same probability as all other alternatives, its entropy will be one hundred percent and one, and it will therefore have no part in the choosing of the option, which appears apparent. This is also described mathematically in broad terms as follows. It will if an index has the same value from the point of view of the m option. Therefore:

$$E_j = -k \sum_{i=1}^{m}[P_i . \ln(P_i)] = -\frac{1}{\ln(m)}[P_1 \ln(P_1), P_2 \ln(P_2),..., P_m \ln(P_m)]$$

$$= -\frac{1}{\ln(m)}\left[\frac{1}{m}\ln\frac{1}{m}, \frac{1}{m}\ln\frac{1}{m},..., \frac{1}{m}\ln\frac{1}{m}\right] = -\frac{1}{\ln(m)}\left[m(\frac{1}{m}\ln\frac{1}{m})\right]$$

$$= -\frac{1}{\ln(m)}\left[1 \times \ln\frac{1}{m}\right] = -\frac{1}{\ln(m)}[-\ln(m)] \Rightarrow E_j = 1$$

That is, such an index is completely entropic, plays no part in option selection, and has no weight, as will be shown.

**Step 4:** The dj value (degree of deviation) is then calculated, which indicates how much usable information the relevant index (dj) offers to the decision maker. The more dissimilar the competing choices are in terms of that index, the closer the measured values of the index are to each other.

$$d_j = 1 - E_j, j = 1,2,...,n \tag{4}$$

**Step 5:** The following equation is used to calculate the weight of each index:

$$W_j = \frac{d_j}{\sum_{j=1}^{m} d_j}, j = 1,2,...,n \tag{5}$$

One of the suitable ways for weighing the criteria is the method presented in this study (WBI). This technique is data-driven, similar to the Shannon entropy method, in that the weight of the criterion is determined by the relationship between the decision matrix data. This approach uses the same reasoning as the entropy method of data dispersion, thus the criteria with the most dispersed data will be more significant for decision making. As a result, he is likely to gain weight.

Low computing burden and simplicity are two advantages of the suggested technique over the Shannon entropy method. Because the data does not need to be normalized. Another benefit of this approach over entropy is that it can be used to calculate negative data, whereas Ln in the Shannon entropy method cannot calculate. Indifference indicates that criteria with zero dispersion has no influence on the decision maker's choice. **DWM** stands for Dispersion-based Weighting Method.

## 2.2. Dispersion-based Weighting Method

**Step 1:** The decision matrix is constructed in the same way as Shannon's entropy technique:

$$X_{ij} = [x_{ij}]_{n \times m} = \begin{bmatrix} x_{11} & x_{12} & \cdots & x_{1m} \\ x_{21} & x_{22} & \cdots & x_{2m} \\ \vdots & \vdots & \vdots & \vdots \\ x_{n1} & x_{n2} & \cdots & x_{nm} \end{bmatrix} \quad (6)$$

**Step 2:** Calculate each criterion's average $\mu_j$, which is determined as follows:

$$\mu_j = \frac{\sum_{i=1}^{n} x_{ij}}{n}, i = 1, 2, \dots, n, j = 1, 2, \dots, m$$

Where n is the number of options and m is the number of criteria.

**Step 3:** Calculate each criterion's standard deviation $s_j$ which is determined as follows:

$$s_j = \sqrt{\frac{(\sum_{i=1}^{n} x_{ij} - \mu_j)^2}{n}}, i = 1, 2, \dots, n, j = 1, 2, \dots, m$$

**Step 4:** Calculate each criterion's coefficient of variation $vc_j$ as follows:

$$vc_j = \frac{s_j}{\mu_j}, j = 1, 2, \dots, m$$

**Step 5:** Use the following equation to calculate the weight of each criterion $w_j$:

$$w_j = \frac{vc_j}{\sum_{j=1}^{m} vc_j}, j = 1, 2, \dots, m$$

## 3- Numerical examples

In this part, we look at two distinct instances to demonstrate the validity of the suggested technique and to show how it differs from the Shannon entropy method.

**Example 3. 1**

Assume a recent university graduate wishes to select one of four occupations based on five indications. Income, social position, hard labor, distance, and social security are all indicators. Table 1 shows the worth of each position in terms of each metric.

**Table1.** Decision Matrix

| Cj<br>Ai | Income<br>C1 | Social image<br>C2 | Hard work<br>C3 | Distance<br>C4 | Security<br>C5 |
|---|---|---|---|---|---|
| **A1** | 15 | High | Relatively high | 10 | High |
| **A2** | 12 | Medium | Medium | 3 | Extremely high |
| **A3** | 20 | Extremely high | High | 30 | Medium |
| **A4** | 30 | Low | Extremely high | 1 | Low |

Two of the five accessible indications (C1, C4) are quantitative, while the others are qualitative, according to this choice matrix. The majority of indicators in MADM models are on various scales and are frequently at odds with one another. Frequently, there is no best choice (ideally for each indicator). Furthermore, certain indications have a good and bad side to them. As a result, under a MADM model, the best choice would be subjective option A, which would offer the best value for each indicator. In the vast majority of situations, getting to A is impossible. In any event, selecting the most suited alternative will be pretty simple.

**3.1.1. Using the Shannon entropy technique to calculate weight:**

**Step 1:** The criterion option matrix is created. The qualitative values are now transformed to numeric values using the Likert scale (Table 2).

**Step 2:** Table 3 shows the normalization and computation of the criterion-option matrix using Equation (2) and computation $p_{ij}$.

**Steps 3-5:** Calculation of all criteria based on entropy. Table 4 combines the results of the entropy calculations, covering Steps 3 through 5. Equations (3) through (6) contain calculating equations (5). The goal values of the financial ratios are shown in step 5 by rj.

**Table2.** Decision matrix

| | Income | Social Image | Hard Work | Distance | Security |
|---|---|---|---|---|---|

| Cj<br>Ai | C1 | C2 | C3 | C4 | C5 |
|---|---|---|---|---|---|
| **A1** | 15 | 6 | 3 | 10 | 6 |
| **A2** | 12 | 4 | 4 | 3 | 7 |
| **A3** | 20 | 7 | 2 | 30 | 4 |
| **A4** | 30 | 2 | 1 | 1 | 2 |

**Table 3.** Normalized matrix

| Cj<br>Ai | Income<br>C1 | Social Image<br>C2 | Hard Work<br>C3 | Distance<br>C4 | Security<br>C5 |
|---|---|---|---|---|---|
| **A1** | 0.1948 | 0.3158 | 0.3000 | 0.2273 | 0.1538 |
| **A2** | 0.1558 | 0.2105 | 0.4000 | 0.0682 | 0.0769 |
| **A3** | 0.2597 | 0.3684 | 0.2000 | 0.6818 | 0.3077 |
| **A4** | 0.3896 | 0.1053 | 0.1000 | 0.0227 | 0.4615 |

### 3.1.2. Calculation of weight using the proposed method (DWM):

**Step 1:** The criterion-option matrix is constructed in the same manner as the Shannon entropy technique.

**Steps 2-5:** Data mean, standard deviation, and coefficient of standard deviation of data were calculated, followed by the weight indicated in Table 5 of the findings.

**Example 3. 2**

The data in this section is based on the Lee article from 2017 [23], which includes examples of four container firms that are in the top 20 in the world: Evergreen Shipping Company, Yang Ming Shipping Company, Hanjin Shipping Company, and Hyundai Merchant Marine Company.

**Table 4**. Entropy calculation for all criteria

| | Step3-lnpij | | | | | Step4 | | Step5 |
|---|---|---|---|---|---|---|---|---|
| Criteria | A1 | A2 | A3 | A4 | $\sum p_{ij} * \ln p_{ij}$ | eij | 1-eij | wj |
| Income | -1.4816 | -1.20397 | -1.15268 | -1.63576 | -1.32575 | 0.9563 | 0.0437 | 0.06325 |
| Social image | -2.68558 | -0.91629 | -1.55814 | -1.8589 | -1.29689 | 0.9355 | 0.0645 | 0.093399 |
| Hard work | -0.38299 | -1.60944 | -0.99853 | -1.34807 | -1.27985 | 0.9232 | 0.0768 | 0.111198 |
| Distance | -3.78419 | -2.30259 | -2.25129 | -0.94261 | -0.86697 | 0.6254 | 0.3746 | 0.542538 |
| Security | -1.4816 | -1.20397 | -1.15268 | -1.63576 | -1.20479 | 0.8691 | 0.1309 | 0.189615 |

**Table 5.** DWM calculation for all criteria

|  | Step1 | | | | Step2 | Step3 | Step4 | Step5 |
|---|---|---|---|---|---|---|---|---|
| Criteria | A1 | A2 | A3 | A4 | $\mu_j$ | $S_j$ | $CV_j$ | $W_j$ |
| Income | 15 | 12 | 20 | 30 | 19.25 | 6.832825 | 0.354952 | 0.124993 |
| Social image | 6 | 4 | 7 | 2 | 4.75 | 1.920286 | 0.404271 | 0.14236 |
| Hard work | 3 | 4 | 2 | 1 | 2.5 | 1.118034 | 0.447214 | 0.157482 |
| Distance | 10 | 3 | 30 | 1 | 11 | 11.46734 | 1.042486 | 0.367101 |
| Security | 6 | 7 | 4 | 2 | 3.25 | 1.920286 | 0.590857 | 0.208065 |

**Table 6.** Comparison of the results of the entropy and proposed method

| Entropy weight | Rank | WBI Weight | Rank |
|---|---|---|---|
| 0.06325 | 5 | 0.124993 | 5 |
| 0.093399 | 4 | 0.14236 | 4 |
| 0.111198 | 3 | 0.157482 | 3 |
| 0.542538 | 1 | 0.367101 | 1 |
| 0.189615 | 2 | 0.208065 | 2 |

### 3.2.1. Shannon entropy technique

**Step 1:** Matrix of Criteria and Options The goal of this phase is to gather the information provided in Table 7. The financial performance of the container transport firms reviewed on their official websites is shown in this table. The last column shows the financial ratio's overall performance, which will be used in the following stage.

**Table 7.** Performance of container shipping companies in financial ratios

|  | Step1-performancE | YM | EG | HMM | HJ | SUM |
|---|---|---|---|---|---|---|
| F1 | Current ratio (%) | 120.30 | 110.99 | 55.24 | 98.26 | 384.79 |
| F2 | Times interest coverage ratio | 4.68 | 1.72 | 0.11 | 1.33 | 7.83 |
| F6 | Gross profit margin | 5.28 | 20.04 | 10.73 | 5.10 | 41.14 |
| F8 | Net profit margin | 10.71 | 11.27 | 1.02 | 7.41 | 30.40 |
| F9 | Income before tax ratio (EBT) (%) | 13.16 | 13.62 | 1.17 | 10.31 | 38.27 |

| F11 | Return on long-term capital | 14.43 | 10.91 | 1.77 | 9.43 | 36.54 |
|---|---|---|---|---|---|---|
| F12 | Return on equity (ROE) | 54.18 | 49.93 | 1.16 | 70.46 | 175.73 |
| F13 | Return on total assets (ROA) | 5.22 | 3.42 | 0.41 | 2.81 | 11.87 |
| F14 | Total Asser Turnover | 1.17 | 0.36 | 0.70 | 0.82 | 3.06 |
| F15 | Fixed Asset Turnover | 2.62 | 1.10 | 1.10 | 1.23 | 6.10 |
| F16 | Stockholder's equity turnover | 1.85 | 0.68 | 7.66 | 5.82 | 16.01 |
| F17 | Total liabilities turnover | 3.16 | 0.78 | 0.78 | 0.96 | 5.67 |
| F18 | Working Capital Turnover | 251.40 | 237.25 | 216.23 | 1.24 | 706.12 |
| F24 | Cash flow to net income ratio | 330.74 | 624.87 | 1.12 | 1035.17 | 1991.90 |
| F25 | Cash flow adequacy ratio | 71.14 | 32.07 | 4.54 | 5.15 | 112.90 |
| F19 | Total debt to asset ratio (Debt ratio) (%) | 36.99 | 46.67 | 90.81 | 85.88 | 260.36 |
| F3 | Equity to total assets ratio (Equity ratio) | 0.63 | 0.53 | 0.09 | 0.14 | 1.40 |
| F5 | Fixed Assets to Long-term Fund Ratio | 0.56 | 0.41 | 1.06 | 0.85 | 2.87 |
| F20 | Debt to Equity Ratio | 0.59 | 0.88 | 988.50 | 6.08 | 996.06 |
| F21 | Long-term debt to Equity ratio | 0.27 | 0.52 | 5.34 | 4.56 | 10.68 |
| F22 | Long-term debt to Long-term Capital ratio | 0.21 | 0.34 | 0.84 | 0.82 | 2.21 |

Notes: EG Evergreen Shipping Company, YM Yang Ming Shipping Company, HJ Hanjin Shipping Company, and HMM Hyundai Merchant Marine Company

**Step 2:** Finding Pij. The total of the normalized financial ratios computed in Table 3 is used to illustrate this phase in Table 3.

**Steps 3-5:** Entropy and total criterion calculation Table 9 combines the results of entropy calculations from Step 3 to Step 5. Calculation equations may be found in Equations (3) through (5). The outcomes of step 5, rj, indicate the financial ratios' goal values.

**Table 8.** Normalized financial ratio performance

| | Step1-performancE | YM | EG | HMM | HJ |
|---|---|---|---|---|---|
| F1 | Current ratio (%) | 0.313 | 0.288 | 0.144 | 0.255 |
| F2 | Times interest coverage ratio | 0.597 | 0.219 | 0.014 | 0.170 |
| F6 | Gross profit margin | 0.128 | 0.487 | 0.261 | 0.124 |
| F8 | Net profit margin | 0.352 | 0.371 | 0.033 | 0.244 |
| F9 | Income before tax ratio (EBT) (%) | 0.344 | 0.356 | 0.031 | 0.270 |
| F11 | Return on long-term capital | 0.395 | 0.299 | 0.048 | 0.258 |
| F12 | Return on equity (ROE) | 0.308 | 0.284 | 0.007 | 0.401 |
| F13 | Return on total assets (ROA) | 0.440 | 0.288 | 0.035 | 0.237 |
| F14 | Total Asser Turnover | 0.382 | 0.119 | 0.230 | 0.269 |
| F15 | Fixed Asset Turnover | 0.429 | 0.181 | 0.187 | 0.202 |
| F16 | Stockholder's equity turnover | 0.116 | 0.042 | 0.478 | 0.363 |
| F17 | Total liabilities turnover | 0.577 | 0.137 | 0.137 | 0.169 |
| F18 | Working Capital Turnover | 0.356 | 0.336 | 0.306 | 0.002 |
| F24 | Cash flow to net income ratio | 0.166 | 0.314 | 0.001 | 0.520 |
| F25 | Cash flow adequacy ratio | 0.630 | 0.284 | 0.040 | 0.046 |
| F19 | Total debt to asset ratio (Debt ratio) (%) | 0.142 | 0.179 | 0.349 | 0.330 |
| F3 | Equity to total assets ratio (Equity ratio) | 0.451 | 0.382 | 0.066 | 0.101 |
| F5 | Fixed Assets to Long-term Fund Ratio | 0.195 | 0.141 | 0.368 | 0.296 |
| F20 | Debt to Equity Ratio | 0.001 | 0.001 | 0.992 | 0.006 |
| F21 | Long-term debt to Equity ratio | 0.025 | 0.048 | 0.500 | 0.427 |
| F22 | Long-term debt to Long-term Capital ratio | 0.095 | 0.154 | 0.380 | 0.370 |

### 3.2.2. Calculation of weight using the proposed method (DWM):

**Step 1:** The criterion-option matrix is constructed in the same manner as the Shannon entropy technique.

**Steps 2-5:** Data mean, standard deviation, and coefficient of standard deviation of data were calculated, followed by the weight indicated in Table 10 of the findings.

**Table 9.** Entropy calculation for all criteria

| | | Steps | | | | Step3-$\ln p_{ij}$ | Step4 | | Step5 | |
|---|---|---|---|---|---|---|---|---|---|---|
| | Financial ratios | YM | EG | HMM | HJ | $\sum p_{ij} * \ln p_i$ | $e_{ij}$ | $1-e_{ij}$ | $w_j$ | |
| F1 | Current ratio (%) | -1.1627 | -1.2433 | -1.9410 | -1.3651 | -1.3494 | 0.9734 | 0.0266 | 0.0065 | |
| F2 | Times interest coverage ratio | -0.5153 | -1.5164 | -4.2940 | -1.7748 | -1.0001 | 0.7214 | 0.2789 | 0.0680 | |
| F6 | Gross profit margin | -2.0540 | -0.7194 | -1.3436 | -2.0887 | -1.2230 | 0.8822 | 0.1178 | 0.0288 | |
| F8 | Net profit margin | -1.0435 | -0.9927 | -3.3973 | -1.4116 | -1.1932 | 0.8607 | 0.1393 | 0.0340 | |
| F9 | Income before tax ratio (EBT) (%) | -1.0673 | -1.0332 | -3.4836 | -1.3111 | -1.1951 | 0.8621 | 0.1379 | 0.0337 | |
| F11 | Return on long-term capital | -0.9293 | -1.2085 | -3.0278 | -1.3544 | -1.2240 | 0.8829 | 0.1171 | 0.0286 | |
| F12 | Return on equity (ROE) | -1.1767 | -1.2582 | -5.0170 | -0.9140 | -1.1200 | 0.8079 | 0.1921 | 0.0469 | |
| F13 | Return on total assets (ROA) | -0.8205 | -1.2438 | -3.3545 | -1.4414 | -1.1780 | 0.8497 | 0.1503 | 0.0367 | |
| F14 | Total Asser Turnover | -0.9614 | -2.1324 | -1.4682 | -1.3139 | -1.3117 | 0.9462 | 0.0538 | 0.0131 | |
| F15 | Fixed Asset Turnover | -0.8452 | -1.7084 | -1.6747 | -1.5994 | -1.3093 | 0.9445 | 0.0555 | 0.0136 | |

| | | | | | | | | | |
|---|---|---|---|---|---|---|---|---|---|
| F16 | Stockholder's equity turnover | -2.1558 | -3.1601 | -0.7372 | -1.0125 | -1.1043 | 0.7966 | 0.2034 | 0.0497 |
| F17 | Total liabilities turnover | -0.5844 | -1.9878 | -1.9893 | -1.7792 | -1.1705 | 0.8443 | 0.1557 | 0.0380 |
| F18 | Working Capital Turnover | -1.0327 | -1.0907 | -1.1834 | -6.3482 | -1.1076 | 0.7990 | 0.2010 | 0.0491 |
| F24 | Cash flow to net income ratio | -1.7955 | -1.1593 | -7.4836 | -0.6545 | -1.0062 | 0.7258 | 0.2742 | 0.0670 |
| F25 | Cash flow adequacy ratio | -0.4618 | -1.2585 | -3.2136 | -3.0880 | -0.9185 | 0.6626 | 0.3374 | 0.0824 |
| F19 | Total debt to asset ratio (Debt ratio) (%) | -1.9514 | -1.7190 | -1.0533 | -1.1091 | -1.3186 | 0.9512 | 0.0488 | 0.0119 |
| F3 | Equity to total assets ratio (Equity ratio) | -0.7958 | -0.9626 | -2.7213 | -2.2917 | -1.1374 | 0.8205 | 0.1795 | 0.0438 |
| F5 | Fixed Assets to Long-term Fund Ratio | -1.6373 | -1.9572 | -0.9988 | -1.2177 | -1.3231 | 0.9544 | 0.0456 | 0.0111 |
| F20 | Debt to Equity Ratio | -7.4364 | -7.0372 | -0.0076 | -5.0982 | -0.0493 | 0.0355 | 0.9645 | 0.2355 |
| F21 | Long-term debt to Equity ratio | -3.6868 | -3.0303 | -0.6933 | -0.8517 | -0.9487 | 0.6844 | 0.3156 | 0.0771 |
| F22 | Long-term debt to Long-term Capital ratio | -2.3502 | -1.8726 | -0.9663 | -0.9930 | -1.2475 | 0.8999 | 0.1001 | 0.0245 |

**Table 10.** DWM calculation for all criteria

| | Steps | Step1 | | | | Step2 | Step3 | Step4 | Step5 |
|---|---|---|---|---|---|---|---|---|---|
| | Financial ratios | YM | EG | HMM | HJ | $\mu_j$ | Sj | Vcj | wj |
| F1 | Current ratio (%) | 120.30 | 110.99 | 55.24 | 98.26 | 96.1975 | 24.90742 | 0.25892 | 0.0193 |
| F2 | Times interest coverage ratio | 4.68 | 1.72 | 0.11 | 1.33 | 1.96 | 1.678943 | 0.856604 | 0.0638 |
| F6 | Gross profit margin | 5.28 | 20.04 | 10.73 | 5.10 | 10.2875 | 6.068201 | 0.589862 | 0.0439 |
| F8 | Net profit margin | 10.71 | 11.27 | 1.02 | 7.41 | 7.6025 | 4.076563 | 0.536214 | 0.0399 |
| F9 | Income before tax ratio (EBT) (%) | 13.16 | 13.62 | 1.17 | 10.31 | 9.565 | 5.009943 | 0.506162 | 0.0390 |
| F11 | Return on long-term capital | 14.43 | 10.91 | 1.77 | 9.43 | 9.135 | 4.623794 | 0.523779 | 0.0377 |
| F12 | Return on equity (ROE) | 54.18 | 49.93 | 1.16 | 70.46 | 43.9325 | 25.85622 | 0.588544 | 0.0438 |
| F13 | Return on total assets (ROA) | 5.22 | 3.42 | 0.41 | 2.81 | 2.965 | 1.720763 | 0.580358 | 0.0432 |

| | | | | | | | | |
|---|---|---|---|---|---|---|---|---|
| F14 | Total Asser Turnover | 1.17 | 0.36 | 0.70 | 0.82 | 0.7625 | 0.289515 | 0.379691 | 0.0283 |
| F15 | Fixed Asset Turnover | 2.62 | 1.10 | 1.14 | 1.23 | 1.5225 | 0.635389 | 0.417332 | 0.0311 |
| F16 | Stockholder's equity turnover | 1.85 | 0.68 | 7.66 | 5.82 | 4.0025 | 2.843979 | 0.710551 | 0.0529 |
| F17 | Total liabilities turnover | 3.16 | 0.78 | 0.78 | 0.96 | 1.42 | 1.007274 | 0.709348 | 0.0528 |
| F18 | Working Capital Turnover | 251.40 | 237.25 | 216.23 | 1.24 | 176.53 | 101.9744 | 0.57766 | 0.0430 |
| F24 | Cash flow to net income ratio | 330.74 | 624.87 | 1.12 | 1035.17 | 497.975 | 380.6288 | 0.764353 | 0.0569 |
| F25 | Cash flow adequacy ratio | 71.14 | 32.07 | 4.54 | 5.15 | 28.225 | 27.15657 | 0.962146 | 0.0716 |
| F19 | Total debt to asset ratio (Debt ratio) (%) | 36.99 | 46.67 | 90.81 | 85.88 | 65.0875 | 23.57249 | 0.362166 | 0.0270 |
| F3 | Equity to total assets ratio (Equity ratio) | 0.63 | 0.53 | 0.09 | 0.14 | 0.3475 | 0.235836 | 0.678666 | 0.0505 |
| F5 | Fixed Assets to Long-term Fund Ratio | 0.56 | 0.41 | 1.06 | 0.85 | 0.72 | 0.252091 | 0.350127 | 0.0261 |
| F20 | Debt to Equity Ratio | 0.59 | 0.88 | 988.50 | 6.08 | 249.0125 | 426.9489 | 1.714568 | 0.1276 |
| F21 | Long-term debt to Equity ratio | 0.27 | 0.52 | 5.34 | 4.56 | 2.6725 | 2.295837 | 0.85906 | 0.0639 |
| F22 | Long-term debt to Long-term Capital ratio | 0.21 | 0.34 | 0.84 | 0.82 | 0.5525 | 0.281369 | 0.509266 | 0.0379 |

**3.3. Comparing results:**

Two numerical examples were provided in this section, and the weights of the criterion were determined using the two weighing methods stated previously (Shannon entropy and the proposed method). Tables 6 and 11 present the findings. All values and rankings of the criterion were the same in the first example with four possibilities, and all values and rankings of the criteria were the same in the second example with 22 options except for 5, which was the same in both ways.

The Pearson correlation coefficient was used to assess the relationship intensity, kind of relationship (direct or inverse), and weight using both techniques to evaluate the effectiveness of the **DWM**. This coefficient ranges from 1 to -1, and it will be equal to zero if there is no link between the two variables. The correlation rate between the values produced by weighting by entropy technique and the suggested method (**DWM**) was calculated using SPSS software and found to be 0.997 in the first case and 0.979 in the second. This demonstrates a high level of connection. The correlation coefficient was calculated using SPSS software, and the findings are shown below.

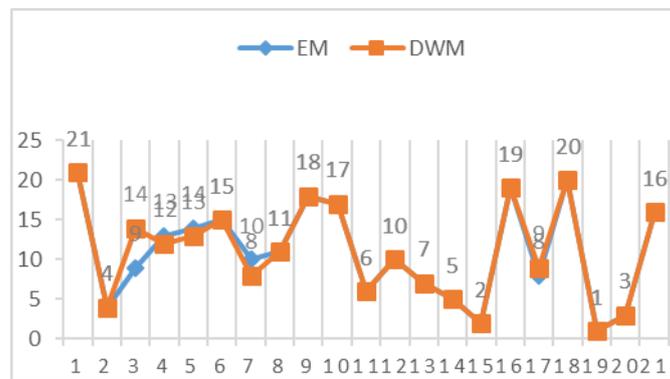

**Fig 1: Comparison of EP and DWM method results**

Table 11. Comparison of the results of the entropy and proposed method

|  | Steps | Entropy Method | | Proposed Method | |
|---|---|---|---|---|---|
|  | Financial ratios | wj | rank | wj | rank |
| F1 | Current ratio (%) | 0.0065 | 21 | 0.0193 | 21 |
| F2 | Times interest coverage ratio | 0.0680 | 4 | 0.0638 | 4 |
| F6 | Gross profit margin | 0.0288 | 14 | 0.0439 | 9 |
| F8 | Net profit margin | 0.0340 | 12 | 0.0399 | 13 |
| F9 | Income before tax ratio (EBT) (%) | 0.0337 | 13 | 0.0390 | 14 |
| F11 | Return on long-term capital | 0.0286 | 15 | 0.0377 | 15 |
| F12 | Return on equity (ROE) | 0.0469 | 8 | 0.0438 | 10 |
| F13 | Return on total assets (ROA) | 0.0367 | 11 | 0.0432 | 11 |
| F14 | Total Asser Turnover | 0.0131 | 18 | 0.0283 | 18 |
| F15 | Fixed Asset Turnover | 0.0136 | 17 | 0.0311 | 17 |
| F16 | Stockholder's equity turnover | 0.0497 | 6 | 0.0529 | 6 |
| F17 | Total liabilities turnover | 0.0380 | 10 | 0.0430 | 10 |
| F18 | Working Capital Turnover | 0.0491 | 7 | 0.0528 | 7 |
| F24 | Cash flow to net income ratio | 0.0670 | 5 | 0.0569 | 5 |
| F25 | Cash flow adequacy ratio | 0.0824 | 2 | 0.0716 | 2 |
| F19 | Total debt to asset ratio (Debt ratio) (%) | 0.0119 | 19 | 0.0270 | 19 |
| F3 | Equity to total assets ratio (Equity ratio) | 0.0438 | 9 | 0.0505 | 8 |
| F5 | Fixed Assets to Long-term Fund Ratio | 0.0111 | 20 | 0.0261 | 20 |
| F20 | Debt to Equity Ratio | 0.2355 | 1 | 0.1276 | 1 |
| F21 | Long-term debt to Equity ratio | 0.0771 | 3 | 0.0639 | 3 |
| F22 | Long-term debt to Long-term Capital ratio | 0.0245 | 16 | 0.0379 | 16 |

**Correlations**

| | | VAR00007 | VAR00009 |
|---|---|---|---|
| VAR00007 | Pearson Correlation | 1 | .997** |
| | Sig. (2-tailed) | | .000 |
| | N | 5 | 5 |
| VAR00009 | Pearson Correlation | .997** | 1 |
| | Sig. (2-tailed) | .000 | |
| | N | 5 | 5 |

**. Correlation is significant at the 0.01 level (2-tailed).

**Correlations**

| | | VAR00002 | VAR00004 |
|---|---|---|---|
| VAR00002 | Pearson Correlation | 1 | .980** |
| | Sig. (2-tailed) | | .000 |
| | N | 21 | 21 |
| VAR00004 | Pearson Correlation | .980** | 1 |
| | Sig. (2-tailed) | .000 | |
| | N | 21 | 21 |

**. Correlation is significant at the 0.01 level (2-tailed).

## 4. Conclusion

In this article we suggest an acceptable technique of weighing by criteria termed Dispersion-based Weighting Method (**DWM**). The weight was calculated in five steps using this approach. To show the applicability of this novel technique, several numerical examples were utilized. We also compared the outcomes of **DWM** with entropy Shannon using a variety of assessment criteria, finding that **DWM** outperformed Shannon entropy. The weight of the criterion is determined by the relationship between the decision matrix data, comparable to the data-driven Shannon entropy approach. This approach uses the same reasoning as the entropy method of data dispersion. As a result, the larger the scatter value of a criterion for several choices, the more weight that criterion is given. **DWM** offers a number of great features that make weighing simple and convenient:

- Low computational burden and simplicity are two benefits of the suggested technique over the Shannon entropy method. Because the data does not need to be normalized.
- Another benefit of this approach over entropy is that it can handle negative data, whereas the Shannon ln entropy method cannot.

Integration of this method with multi-criteria decision-making systems like Topsis might be an intriguing future research project. This method presented is in crisp space and can be developed in future research based on fuzzy logic. We also suggest that you apply the suggested approach to some additional real-world applications and compare the findings to those of other MCDM methods in order to enhance the method's validity and verify its applicability and utility.